%

\magnification=1200 \vsize=18cm \voffset=1cm 
\hoffset=-.1cm



\hsize=11.25cm
\parskip 0pt
\parindent=12pt

\catcode'32=9

\font\tenpc=cmcsc10
\font\eightpc=cmcsc8
\font\eightrm=cmr8
\font\eighti=cmmi8
\font\eightsy=cmsy8
\font\eightbf=cmbx8
\font\eighttt=cmtt8
\font\eightit=cmti8
\font\eightsl=cmsl8
\font\sixrm=cmr6
\font\sixi=cmmi6
\font\sixsy=cmsy6
\font\sixbf=cmbx6

\skewchar\eighti='177 \skewchar\sixi='177
\skewchar\eightsy='60 \skewchar\sixsy='60

\catcode`@=11

\def\tenpoint{%
 \textfont0=\tenrm \scriptfont0=\sevenrm \scriptscriptfont0=\fiverm
 \def\rm{\fam\z@\tenrm}%
 \textfont1=\teni \scriptfont1=\seveni \scriptscriptfont1=\fivei
 \def\oldstyle{\fam\@ne\teni}%
 \textfont2=\tensy \scriptfont2=\sevensy \scriptscriptfont2=\fivesy
 \textfont\itfam=\tenit
 \def\it{\fam\itfam\tenit}%
 \textfont\slfam=\tensl
 \def\sl{\fam\slfam\tensl}%
 \textfont\bffam=\tenbf \scriptfont\bffam=\sevenbf
 \scriptscriptfont\bffam=\fivebf
 \def\bf{\fam\bffam\tenbf}%
 \textfont\ttfam=\tentt
 \def\tt{\fam\ttfam\tentt}%
 \abovedisplayskip=12pt plus 3pt minus 9pt
 \abovedisplayshortskip=0pt plus 3pt
 \belowdisplayskip=12pt plus 3pt minus 9pt
 \belowdisplayshortskip=7pt plus 3pt minus 4pt
 \smallskipamount=3pt plus 1pt minus 1pt
 \medskipamount=6pt plus 2pt minus 2pt
 \bigskipamount=12pt plus 4pt minus 4pt
 \normalbaselineskip=12pt
 \setbox\strutbox=\hbox{\vrule height8.5pt depth3.5pt width0pt}%
 \let\bigf@ntpc=\tenrm \let\smallf@ntpc=\sevenrm
 \let\petcap=\tenpc
 \normalbaselines\rm}

\def\eightpoint{%
 \textfont0=\eightrm \scriptfont0=\sixrm \scriptscriptfont0=\fiverm
 \def\rm{\fam\z@\eightrm}%
 \textfont1=\eighti \scriptfont1=\sixi \scriptscriptfont1=\fivei
 \def\oldstyle{\fam\@ne\eighti}%
 \textfont2=\eightsy \scriptfont2=\sixsy \scriptscriptfont2=\fivesy
 \textfont\itfam=\eightit
 \def\it{\fam\itfam\eightit}%
 \textfont\slfam=\eightsl
 \def\sl{\fam\slfam\eightsl}%
 \textfont\bffam=\eightbf \scriptfont\bffam=\sixbf
 \scriptscriptfont\bffam=\fivebf
 \def\bf{\fam\bffam\eightbf}%
 \textfont\ttfam=\eighttt
 \def\tt{\fam\ttfam\eighttt}%
 \abovedisplayskip=9pt plus 2pt minus 6pt
 \abovedisplayshortskip=0pt plus 2pt
 \belowdisplayskip=9pt plus 2pt minus 6pt
 \belowdisplayshortskip=5pt plus 2pt minus 3pt
 \smallskipamount=2pt plus 1pt minus 1pt
 \medskipamount=4pt plus 2pt minus 1pt
 \bigskipamount=9pt plus 3pt minus 3pt
 \normalbaselineskip=9pt
 \setbox\strutbox=\hbox{\vrule height7pt depth2pt width0pt}%
 \let\bigf@ntpc=\eightrm \let\smallf@ntpc=\sixrm
 \let\petcap=\eightpc
 \normalbaselines\rm}
\catcode`@=12

\newif\ifpagetitre
\newtoks\auteurcourant \auteurcourant={\hfil}
\newtoks\titrecourant \titrecourant={\hfil}


\newif\ifpagetitre
\newtoks\auteurcourant \auteurcourant={\hfil}
\newtoks\titrecourant \titrecourant={\hfil}

\def\appeln@te{}
\def\vfootnote#1{\def\@parameter{#1}\insert\footins\bgroup\eightpoint
 \interlinepenalty\interfootnotelinepenalty
 \splittopskip\ht\strutbox 
 \splitmaxdepth\dp\strutbox \floatingpenalty\@MM
 \leftskip\z@skip \rightskip\z@skip
 \ifx\appeln@te\@parameter\indent \else{\noindent #1\ }\fi
 \footstrut\futurelet\next\fo@t}

\pretolerance=500 \tolerance=1000 \brokenpenalty=5000
\newdimen\hmargehaute \hmargehaute=0cm
\newdimen\lpage \lpage=13.3cm
\newdimen\hpage \hpage=20cm
\newdimen\lmargeext \lmargeext=1cm
\hsize=11.25cm
\vsize=18cm
\parskip 0pt
\parindent=12pt

\def\margehaute{\vbox to \hmargehaute{\vss}}%
\def\margebasse{\vss}

\output{\shipout\vbox to \hpage{\margehaute\nointerlineskip
 \corpsdepage\margebasse}
 \advancepageno \global\pagetitrefalse
 \ifnum\outputpenalty>-20000 \else\dosupereject\fi}

\def\corpsdepage{\hbox to \lpage{\hss\pagetexte\hskip\lmargeext}}
\def\pagetexte{\vbox{\makeheadline\pagebody\makefootline}}
\headline={\ifpagetitre\titleheadline \else
 \ifodd\pageno\rightheadline \else\leftheadline\fi\fi}
\def\leftheadline{\eightpoint\hfil\the\auteurcourant\hfil}
\def\rightheadline{\eightpoint\hfil\the\titrecourant\hfil}
\def\titleheadline{\hfill}
\pagetitretrue

\def\footnoterule{\kern-6\p@
 \hrule width 2truein \kern 5.6\p@} 

\tenpoint

\font\tengoth=eufm10
\def\goth#1{\hbox{\tengoth #1}}

\catcode`\@=11
\def\pc#1#2|{{\bigf@ntpc #1\penalty \@MM\hskip\z@skip\smallf@ntpc%
   \uppercase{#2}}}
\catcode`\@=12

\def\pointir{\discretionary{.}{}{.\kern.35em---\kern.7em}\nobreak
  \hskip 0em plus .3em minus .4em }

\def\qed{\quad\raise -2pt\hbox{\vrule\vbox to 10pt{\hrule width 4pt
  \vfill\hrule}\vrule}}

\def\cqfd{\penalty 500 \hbox{\qed}\par\smallskip}

\def\rem#1|{\par\medskip{{\it #1}.\quad}}

\def\vspace[#1]{\noalign{\vskip#1}}

\def\abstract#1{\vbox{\eightpoint\narrower\narrower
\pc ABSTRACT|.\quad #1}}

\def\section#1{\goodbreak\par\vskip .3cm\centerline{\bf #1}
  \par\nobreak\vskip 3pt }

\long\def\th#1|#2\endth{\par\medbreak
  {\petcap #1\pointir}{\it #2}\par\medbreak}

\def\article#1|#2|#3|#4|#5|#6|#7|
   {{\leftskip=7mm\noindent
    \hangindent=2mm\hangafter=1
    \llap{{\tt [#1]}\hskip.35em}{#2}.\quad
    #3, {\sl #4}, {\bf #5} ({\oldstyle #6}),
    pp.\nobreak\ #7.\par}}

\def\livre#1|#2|#3|#4|
   {{\leftskip=7mm\noindent
   \hangindent=2mm\hangafter=1
   \llap{{\tt [#1]}\hskip.35em}{#2}.\quad
   {\sl #3}, #4.\par}}

\def\divers#1|#2|#3|
   {{\leftskip=7mm\noindent
   \hangindent=2mm\hangafter=1
    \llap{{\tt [#1]}\hskip.35em}{#2}.\quad
    #3.\par}}

\def\proof{\par{\it Proof}.\quad}
\def\qed{\quad\raise -2pt\hbox{\vrule\vbox to 10pt{\hrule width 4pt
\vfill\hrule}\vrule}}

\def\cqfd{\penalty 500 \hbox{\qed}\par\smallskip}

\font\tengoth=eufm10

\def\des{\mathop{\rm des}\nolimits}
\def\exc{\mathop{\rm exc}\nolimits}
\def\dez{\mathop{\rm dez}}
\def\maj{\mathop{\rm maj}\nolimits}
\def\maz{\mathop{\rm maz}}
\def\maf{\mathop{\rm maf}\nolimits}
\def\mafz{\mathop{\rm mafz}\nolimits}
\def\imaj{\mathop{\rm imaj}\nolimits}
\def\pix{\mathop{\rm pix}\nolimits}
\def\fix{\mathop{\rm fix}\nolimits}
\def\FIX{\mathop{\rm FIX}\nolimits}

\def\Der{\mathop{\rm Der}\nolimits}
\def\FIX{\mathop{\hbox{\eightrm FIX}}}

\def\DES{\mathop{\hbox{\eightrm DES}}}
\def\DE2{\mathop{\hbox{\eightrm DEZ}}}

\def\Desar{\mathop{\rm Desar}\nolimits}
\def\Idesar{\mathop{\rm Idesar}\nolimits}
\def\mag{\mathop{\rm mag}\nolimits}
\def\PIX{\mathop{\hbox{\eightrm PIX}}}
\def\DEZ{\mathop{\hbox{\eightrm DEZ}}}
\def\RISE{\mathop{\hbox{\eightrm RISE}}\nolimits}
\def\Rise{\mathop{\hbox{\eightrm RISE}}\nolimits}
\def\bfPhi{{\bf \Phi}}
\def\bfPsi{{\bf \Psi}}
\def\rmF{{\rm F}}
\def\Sh{\mathop{\rm Sh}\nolimits}
\def\Zero{\mathop{\rm Zero}\nolimits}
\def\zero{\mathop{\rm zero}\nolimits}

\def\red{\mathop{\rm red}\nolimits}
\def\ZDer{\mathop{\rm Z\hskip-.8pt Der}\nolimits}

\def\Longmapsto#1{\hbox to #1{$\mapstochar$\rightarrowfill}}

\def\RISE{\mathop{\hbox{\eightrm RISE}}\nolimits}
\def\RIZE{\mathop{\hbox{\eightrm RIZE}}\nolimits}
\def\DES{\mathop{\hbox{\eightrm DES}}\nolimits}
\def\Pos{\mathop{\rm Pos}\nolimits}

%
\catcode`\@=11

\def\matrice#1{\null \,\vcenter {\normalbaselines \m@th
\ialign {\hfil $##$\hfil &&\  \hfil $##$\hfil\crcr
\mathstrut \crcr \noalign {\kern -\baselineskip } #1\crcr
\mathstrut \crcr \noalign {\kern -\baselineskip }}}\,}

\def\petitematrice#1{\left(\null\vcenter {\normalbaselines \m@th
\ialign {\hfil $##$\hfil &&\thinspace  \hfil $##$\hfil\crcr
\mathstrut \crcr \noalign {\kern -\baselineskip } #1\crcr
\mathstrut \crcr \noalign {\kern -\baselineskip }}}\right)}

\catcode`\@=12

\frenchspacing



\titrecourant={TWO TRANSFORMATIONS}
\auteurcourant={DOMINIQUE FOATA AND GUO-NIU HAN}
\rightline{2006/11/29/14:15}
\bigskip
\bigskip
\bigskip

\centerline{\bf Fix-Mahonian Calculus, I: two transformations}
\bigskip
\centerline{
\sl Dominique Foata and Guo-Niu Han}

\bigskip

\bigskip
\hbox{\hskip5cm\vbox{\eightpoint
\eightsl\vbox{\halign{#\hfill\cr
Tu es Petrus, et super hanc petram,\cr
aedificavisti Lacim Uqam tuam.\cr
\noalign{\medskip}
To Pierre Leroux, Montreal, Sept. 2006,\cr
on the occasion of the LerouxFest.\cr}}}}

\bigskip\medskip
\abstract{We construct two bijections of the symmetric group
${\goth S}_n$ onto itself that enable us to show that  three new
three-variable statistics are equidistributed  with classical
statistics involving the number of fixed points. The first one is
equidistributed with the triplet $(\fix,\des,\maj)$, the last two with
$(\fix,\exc,\maj)$, where ``fix," ``des," ``exc" and ``maj" denote the number of
fixed points, the number of descents, the number of excedances and the major
index, respectively.}

\bigskip

\centerline{\bf 1. Introduction} 
\medskip
In this paper {\it Fix-Mahonian Calculus} is understood to
mean the study of multivariable statistics on the symmetric
group~${\goth S}_n$, which involve the number of fixed points ``fix"
as a marginal component. As for the two transformations mentioned
in the title, they make it possible to show that the new statistics
defined below are equidistributed with the classical ones. Those
transformations will be described not directly on~${\goth S}_n$, but
on classes of {\it shuff\kern0.5pt les}, as now introduced.

Let $0\le m\le n$ and let $v$ be a nonempty word of length $m$,
whose letters are {\it positive} integers (with possible repetitions).
Designate by $\Sh(0^{n-m}v)$ the set of all {\it shuff\kern0.5pt les}
of the words $0^{n-m}$ and $v$, that is, the set of all
rearrangements of the juxtaposition product $0^{n-m}v$, whose
longest {\it subword}  of positive letters is~$v$.  Let
$w=x_1x_2\cdots x_n$ be a word from
$\Sh(0^{n-m}v)$. It is convenient to write:
$\Pos w:=v$, $\Zero w:=\{i:1\le i\le n,\,x_i=0\}$, 
$\zero w:=\#\Zero w\ (=n-m)$, so that~$w$ is completely
characterized by the pair $(\Zero w,\Pos w)$.

The {\it descent set}, $\DES w$, and {\it rise set},
$\RISE w$, of~$w$ are respectively defined as being the {\it
subsets}:
$$
\leqalignno{\noalign{\vskip-1pt}
\DES w&:=\{i:1\le i\le n-1, x_i>x_{i+1}\};&(1.1)\cr
\RISE w&:=\{i:1\le i\le n, x_i\le x_{i+1}\}.&(1.2)\cr}
$$

\goodbreak
\noindent
By convention, $x_0=x_{n+1}=+\infty$. The {\it
major index} of~$w$ is defined by
$$
\leqalignno{
\maj w&:=\sum_{i\ge 1} i\quad(i\in \DES w),&(1.3)\cr
\noalign{\hbox{and a new integral-valued statistic ``mafz" by}}
\mafz w&:=\sum_{i\in \Zero w} i
-\sum_{i=1}^{\zero w}i+\maj \Pos w.&(1.4)\cr}
$$

Note that the first three definitions are also valid for each arbitrary
word with nonnegative letters. The link of ``mafz" with the statistic
``maf" introduced in {\tt[CHZ97]} for permutations will be further
mentioned. 

\goodbreak
We shall also be interested in shuffle classes $\Sh(0^{n-m}v)$ when
the word~$v$ is a {\it derangement} of the set
$[\,m\,]:=\{1,2,\ldots,m\}$, that is, when the word
$v=y_1y_2\cdots y_m$ is a permutation of $12\cdots m$ and
$y_i\not=i$ for all~$i$. For short, $v$ is a {\it derangement
of order}~$m$. Let $w=x_1x_2\cdots x_n$ a be word from the
shuffle class $\Sh(0^{n-m}v)$. Then $v=y_1y_2\cdots
y_m=x_{j_1}x_{j_2}\cdots x_{j_m}$ for a certain sequence $1\le
j_1<j_2<\cdots <j_m\le n$. Let ``red" be the increasing bijection of
$\{j_1,j_2,\ldots,j_m\}$ onto $[\,m\,]$. Say that each positive
letter~$x_k$ of~$w$ is {\it excedent} (resp. {\it subexcedent\/})
if and only if
$x_k>\red k$ (resp. $x_k<\red k$). Another kind of rise set, denoted
by $\Rise^\bullet w$, can then be introduced as follows.

Say that $i\in \Rise^\bullet w$ if and only if $1\le i\le n$ and if
one of the following conditions holds (assuming that
$x_{n+1}=+\infty$):

(1) $0<x_i<x_{i+1}$;

(2) $x_i=x_{i+1}=0$;

(3) $x_i=0$ and $x_{i+1}$ is excedent;

(4) $x_i$ is subexcedent and $x_{i+1}=0$.

\noindent
Note that if $x_i=0$ and $x_{i+1}$ is subexcedent, then $i\in
\Rise w\setminus\Rise^\bullet w$, while if $x_i$ is subexcedent and
$x_{i+1}=0$, then $i\in \Rise^\bullet w\setminus\Rise w$.

\medskip
{\it Example}.\quad
Let $v=5\,1\,2\,3\,6\,4$ be a derangement of order~6. Its excedent
letters are 5, 6. Let
$w=5\,0\,1\,2\, 0\,0\,3\,6\,4\in \Sh(0^3v)$. Then,
$\RISE w= \{2,3,5,6,7,9\}$ and
$\RISE^\bullet w=\{3,4 ,5,7,9\}$. 
Also $\mafz w=(2+5+6)-(1+2+3)+\maj (512354)=7+6=13$.

\proclaim Theorem 1.1. For each derangement~$v$ of order~$m$ and
each integer $n\ge m$ the transformation $\bfPhi$ constructed in
Section~$2$ is a bijection of $\Sh(0^{n-m}v)$ onto itself having
the property that
$$
\RISE w=\RISE^\bullet \bfPhi(w)\leqno(1.5)
$$
holds for every $w\in \Sh(0^{n-m}v)$.

\proclaim Theorem 1.2. For each arbitrary word~$v$ of length~$m$
with positive letters and each integer $n\ge m$ the
transformation~${\bf F}_3$ constructed in Section~$4$ is 
a bijection of $\Sh(0^{n-m}v)$ onto itself having
the property that
$$\leqalignno{
\maj w&=\mafz {\bf F}_3(w);&(1.6)\cr
L\,w&= L\,{\bf F}_3(w)\quad (``L" \hbox{ for ``last" or rightmost
letter}); &(1.7)\cr
}
$$
hold for every $w\in \Sh(0^{n-m}v)$.

We emphasize the fact that Theorem 1.1 is restricted to the case
where~$v$ is a derangement, while Theorem 1.2 holds for an
arbitrary word~$v$ with possible repetitions. In Fig.~1 
we can see that ``$\RISE$" and ``$\RISE^\bullet$" 
(resp. ``$\maj$" and ``$\mafz$") are equidistributed on the
shuffle class $\Sh(0^2312)$ (resp. $\Sh(0^2121)$). 

\goodbreak
{
\eightpoint
$$\hbox{\vbox{\offinterlineskip
\halign{\vrule#\hfil\strut&\ \hfil$#$\hfil\ \vrule&\ \hfil$#$\hfil\
\vrule &\ \hfil$#$\hfil\ \vrule&\ \hfil$#\,$\hfil\ \vrule\cr
\noalign{\hrule}
&\hbox{\sevenrm RISE}\, w&w&{\bf \Phi}(w)&\hbox{\sevenrm
RISE}^\bullet{\bf\Phi}(w)\cr
\noalign{\hrule}
&1,2,4,5&0\,0\,3\,1\,2&0\,0\,3\,1\,2&1,2,4,5\cr
\noalign{\hrule}
&1,3,4,5&0\,3\,0\,1\,2&0\,3\,1\,2\,0&1,3,4,5\cr
&&0\,3\,1\,0\,2&0\,3\,0\,1\,2&\cr
\noalign{\hrule}
&1,3,5&0\,3\,1\,2\,0&0\,3\,1\,0\,2&1,3,5\cr
\noalign{\hrule}
&2,3,4,5&3\,0\,0\,1\,2&3\,1\,2\,0\,0&2,3,4,5\cr
\noalign{\hrule}
&2,4,5&3\,0\,1\,0\,2&3\,0\,0\,1\,2&2,4,5\cr
&&3\,1\,2\,0\,0&3\,1\,0\,2\,0&\cr
\noalign{\hrule}
&2,3,5&3\,0\,1\,2\,0&3\,0\,1\,0\,2&2,3,5\cr
\noalign{\hrule}
&3,4,5&3\,1\,0\,0\,2&3\,0\,1\,2\,0&3,4,5\cr
\noalign{\hrule}
&3,5&3\,1\,0\,2\,0&3\,1\,0\,0\,2&3,5\cr
\noalign{\hrule}
\noalign{\medskip}
\multispan5\hfil$\Sh(0^2312)$\hfil\cr
}}\qquad
\vbox{\offinterlineskip\halign{
\vrule#\hfil\strut\ &\ \hfil$#$\hfil\ \vrule&\ \hfil$#$\hfil\
\vrule &\ \hfil$#$\hfil\ \vrule&\ \hfil$#\,$\hfil\ \vrule\cr
\noalign{\hrule}
&\maj w&w&{\bf F}_3(w)&\mafz{\bf F}_3(w)\cr
\noalign{\hrule}
&2&1\,2\,0\,0\,1&0\,0\,1\,2\,1&2\cr
\noalign{\hrule}
&3&0\,1\,2\,0\,1&0\,1\,0\,2\,1&3\cr
\noalign{\hrule}
&4&0\,0\,1\,2\,1&1\,0\,0\,2\,1\,&4\cr
&&1\,0\,2\,0\,1&0\,1\,2\,0\,1&\cr
\noalign{\hrule}
&5&1\,0\,0\,2\,1&1\,0\,2\,0\,1&5\cr
&&1\,2\,1\,0\,0&0\,1\,2\,1\,0&\cr
\noalign{\hrule}
&6&0\,1\,0\,2\,1&1\,2\,0\,0\,1&6\cr
&&1\,2\,0\,1\,0&1\,0\,2\,1\,0&\cr
\noalign{\hrule}
&7&0\,1\,2\,1\,0&1\,2\,0\,1\,0&7\cr
\noalign{\hrule}
&8&1\,0\,2\,1\,0&1\,2\,1\,0\,0&8\cr
\noalign{\hrule}
\noalign{\medskip}
\multispan5\hfil$\Sh(0^2121)$\hfil\cr
}}}
$$

}

\vskip-.4cm
\centerline{Fig. 1}

\medskip
Those two transformations are fully exploited once we know how to
map those shuffle classes onto the symmetric groups. The
permutations from the symmetric group~${\goth S}_n$ will be
regarded as linear words
$\sigma=\sigma(1)\sigma(2)\cdots \sigma(n)$. If
$\sigma$ is such a permutation, let $\FIX\sigma$ denote the set
of its fixed points, i.e., $\FIX\sigma:=\{i:1\le i\le
n,\,\sigma(i)=i\}$ and let $\fix\sigma:=\#\FIX\sigma$. Let
$(j_1,j_2,\ldots, j_{m}$) be the increasing sequence of the
integers~$k$ such that $1\le k\le n$ and
$\sigma(k)\not=k$ and ``red" be the increasing bijection of
$\{j_1,j_2,\ldots,j_{m}\}$ onto $[m]$. The word
$w=x_1x_2\cdots x_n$ derived from $\sigma
=\sigma(1)\sigma(2)\cdots\sigma(n)$ by replacing each fixed point
by~0 and each
other letter $\sigma(j_{k})$ by $\red\sigma(j_{k})$ will be denoted
by~$\ZDer(\sigma)$. Also let
$$
\Der\sigma:=\red\sigma(j_1)\,\red\sigma(j_2)\,\cdots\,\red
\sigma(j_m),\leqno(1.8)
$$
so that $\Der\sigma$ is the word derived from~$\ZDer(\sigma)$
by deleting all the zeros. Accordingly, $\Der\sigma=\Pos\,
\ZDer(\sigma)$. 

It is important to notice that $\Der\sigma$ is a {\it derangement}
of order~$m$. Also $\sigma(j_{k})$ is excedent
in~$\sigma$ ({\it i.e.} $\sigma(j_k)>j_k$) if and only
$\red\sigma(j_k)$ is excedent in $\Der\sigma$
({\it i.e.} $\red\sigma(j_k)>\red j_k$)

Recall that the statistics ``$\DES$," ``$\RISE$" and ``$\maj$" are
also valid for permutations
$\sigma=\sigma(1)\sigma(2)\cdots\sigma(n)$ and that the statistics
``des" ({\it number of descents}) and ``exc" ({\it number of
excedances}) are defined by
$$
\leqalignno{\noalign{\vskip-3pt}
\des\sigma&:=\#\DES\sigma;&(1.9)\cr
\exc\sigma&:=\#\{i:1\le i\le n-1,\,\sigma(i)>i\}.&(1.10)\cr
\noalign{\hbox{We further define:}}
\noalign{\vskip-3pt}
\DEZ \sigma&:=\DES\,  \ZDer(\sigma);
&(1.11)\cr
\RIZE \sigma&:=\RISE\, \ZDer(\sigma);&(1.12)\cr
\dez \sigma&:=\#\DEZ\sigma=\des\, \ZDer(\sigma);&(1.13)\cr
\maz \sigma&:=\maj\, \ZDer(\sigma);&(1.14)\cr
\maf \sigma&:=\mafz\,  \ZDer(\sigma).&(1.15)\cr}
$$

\noindent
As the zeros of $\ZDer(\sigma)$ correspond to the
fixed points of $\sigma$, we also have
$$
\maf\sigma:=\sum_{i\in
\hbox{\sevenrm FIX}\,\sigma}i-\sum_{i=1}^{\fix\sigma}i
+\maj\Der\sigma.\leqno(1.16)
$$

{\it Example}.\ Let
$\sigma=
8\,{\bf2}\,1\,3\,{\bf5}\,{\bf6}\,4\,9\,7$; then
$\DES\sigma=\{1,2,6,8\}$, $\des\sigma=4$, $\maj\sigma=17$,
$\exc\sigma=2$. Furthermore,
$\ZDer(\sigma)=w=5\,0\,1\,2\,0\,0\,3\,6\,4$
and
$\Pos w=\Der\sigma=5\,1\,2\,3\,6\,4$ is a derangement of order~6.
We have $\FIX\sigma=\{2,5,6\}$, $\fix\sigma=3$,
$\DEZ\sigma=\{1,4,8\}$, $\RIZE w=\{2,3,5,6,7,9\}$,
$\dez=3$,
$\maz\sigma=13$ and
$\maf\sigma=(2+5+6)-(1+2+3)+\maj(512364)=7+6=13$.

\medskip
For each $n\ge 0$ let $D_n$ be the set of all derangements
of order~$n$ and~${\goth S}_n^{\Der}$ be the union:\quad
${\goth S}_n^{\Der}:=\bigcup\limits_{m,v}
\Sh(0^{n-m}v)$\quad $(0\le m\le n,\, v\in D_m)$. 

\proclaim Proposition 1.3. The map $\ZDer$ is a bijection of
${\goth S}_n$ onto ${\goth S}_n^{\Der}$ having the following
properties:
$$
\RIZE\sigma=\RISE\, \ZDer(\sigma)\quad{\sl and}\quad
\RISE\sigma=\RISE^\bullet\, \ZDer(\sigma).\leqno(1.17)
$$

\proof
It is evident to verify that $\ZDer$ is bijective and to define
its inverse
$\ZDer^{-1}$. On the other hand, we have $\RIZE=\RISE
\ZDer$
by definition. Finally, let
$w=x_1x_2\cdots x_n=\ZDer(\sigma)$ and
$\sigma(i)<\sigma(i+1)$ for
$1\le i\le n-1$. Four cases are to be considered:

(1) both $i$ and $i+1$ are not fixed by~$\sigma$ and
$0<x_i<x_{i+1}$;

(2) both $i$ and $i+1$ are fixed points and $x_i=x_{i+1}=0$;

(3) $\sigma(i)=i$ and $\sigma(i+1)$ is excedent; then
$x_i=0$ and $x_{i+1}$ is also excedent;

(4) $\sigma(i)<i<i+1=\sigma(i+1)$; then $x_i$ is subexcedent
and $x_{i+1}=0$.

\noindent
We recover the four cases considered in the definition of
$\RISE^\bullet$. The case $i=n$ is banal to study.\cqfd

\goodbreak
\medskip
We next form the two chains:
$$
\leqalignno{
\Phi&:\sigma\buildrel \ZDer\over \mapsto w\buildrel \bfPhi\over
\mapsto w'\ \buildrel \ZDer^{-1}\over \mapsto\ 
\sigma';&(1.18)\cr 
{\rm F}_3&:\sigma\buildrel \ZDer\over \mapsto
w\buildrel {\bf F}_3\over
\mapsto w''\ \buildrel \ZDer^{-1}\over \mapsto\ 
\sigma''.&(1.19)\cr 
}
$$
The next theorem is then a consequence of Theorems 1.1
and~1.2 and Propositions~1.3.

\proclaim Theorem 1.4. The mappings $\Phi$, $\rmF_3$
defined by $(1.18)$ and $(1.19)$ are bijections of ${\goth S}_n$
onto itself and have the following properties
$$\leqalignno{
(\fix, \RIZE,
\Der)\,\sigma&=(\fix,\RISE,\Der)\,\Phi(\sigma);&(1.20)\cr
(\fix, \maz,
\Der)\,\sigma&=(\fix,\maf,\Der)\,{\rm F}_3(\sigma);&(1.21)\cr
(\fix,\maj,\Der)\,\sigma&=(\fix,\maf,\Der)\,{\rm F}_3\circ \Phi^{-1}
(\sigma);&(1.22)
\cr
}
$$
for every $\sigma$ from ${\goth S}_n$.

\goodbreak
\medskip

It is evident that if $\Der\sigma=\Der\tau$ holds for a pair of
permutations $\sigma$, $\tau$ of order~$n$, then
$\exc\sigma=\exc\tau$. Since $\DES\sigma=[\,n\,]\setminus \RISE
\sigma$ and $\DEZ\sigma=[\,n\,]\setminus \RIZE\sigma$ it follows
from (1.20) that
$$
\leqalignno{
(\fix,\DEZ  ,\exc)\,\sigma
&=(\fix,\DES,\exc)\,\Phi(\sigma);&(1.23)\cr
(\fix,\dez ,\maz ,\exc)\,\sigma
&=(\fix,\des,\maj,\exc)\,\Phi(\sigma).&(1.24)\cr
}
$$
On the other hand, (1.23) implies that
$$
(\fix,\maz,\exc)\,\sigma=(\fix,\maf,\exc)\, {\rm F}_3(\sigma).
\leqno(1.25)
$$

As a consequence we obtain the following Corollary.

\proclaim Corollary 1.5. The two triplets
$(\fix,\dez,\maz)$ and $(\fix,\des,\maj)$ are equidistributed
over~${\goth S}_n$. Moreover, the three triplets
$(\fix,\exc,\maz)$, $(\fix,\exc,\maj)$ and $(\fix,\exc,\maf)$
are also equidistributed over~${\goth S}_n$.

The distributions of $(\fix,\des,\maj)$ and
$(\fix,\exc,\maj)$ have been calculated 
by Gessel-Reutenauer ({\tt[GeRe93]}, Theorem~8.4) and by Shareshian
and Wachs {\tt[ShWa06]}, respectively, using the algebra of the
$q$-series (see, e.g., Gasper and Rahman ({\tt[GaRa90]}, chap.~1).
Let
$$
\displaylines{
A_n^{\fix,\des,\maj}(Y,t,q):=\sum_{\sigma\in {\goth S}_n}
Y^{\fix\sigma}t^{\des\sigma}q^{\maj\sigma}\quad(n\ge0);\cr
A_n^{\fix,\exc,\maj}(Y,s,q):=\sum_{\sigma\in {\goth S}_n}
Y^{\fix\sigma}s^{\exc\sigma}q^{\maj\sigma}\quad(n\ge0).\cr}
$$

\goodbreak
Then, they respetively derived the identities:
$$
\displaylines{
(1.26)\ \sum_{n\ge 0}A_n^{\fix,\des,\maj}(Y,t,q){u^n\over
(t;q)_{n+1}} =\sum_{r\ge 0}t^r\Bigl(
1-u\sum_{i=0}^r q^i\Bigr)^{-1} {(u;q)_{r+1}\over
(uY;q)_{r+1}};\hfill\cr
(1.27)\quad \sum_{n\ge 0}
A_n^{\fix,\exc,\maj}(Y,s,q){u^n\over
(q;q)_n} ={(1-sq)e_q(Yu)\over e_q(squ)-sqe_q(u)}.\hfill\cr}
$$
We then know the distributions over~${\goth S}_n$ of the
triplets $(\fix,\dez,\maz)$, $(\fix,\exc,\maz)$ and
$(\fix,\exc,\maf)$. Note that the statistic ``maf" was introduced
by Clarke {\it et al.} {\tt[CHZ97]}. Although it was not explicitly
stated, their bijection ``{\eightrm CHZ}"
of~${\goth S}_n$ onto itself satisfies identity (1.22) when
${\rm F}_3\circ \Phi^{-1}$ is replaced by
``{\eightrm CHZ}."

As is shown in Section 2, the transformation~$\bfPhi$ is described
as a composition product of bijections~$\phi_l$. 
The image $\phi_l(w)$ of each word~$w$ from a shuffle class
$\Sh(0^{n-m}v)$ is obtained by moving its $l$-th zero,
to the right or to the left, depending on its preceding and following
letters. The description of the inverse bijection~$\bfPsi$
of~$\bfPhi$ follows an analogous pattern. The verification of
identity (1.6) requires some attention and is made in Section~3.
The construction of the transformation
${\bf F}_3$ is given in Section~4. Recall that ${\bf F}_3$
maps each shuffle class $\Sh(0^{n-m}v)$ onto itself, the word~$v$
being an {\it arbitrary} word with nonnegative letters. Very much
like the {\it second fundamental
transformation} (see, e.g., [Lo83], p.~201, Algorithm 10.6.1) the
construction of ${\bf F}_3$ is defined by induction on the length of
the words and preserves the {\it rightmost} letter. 

\bigskip
\centerline{\bf 2. The bijection $\bfPhi$} 

\medskip
Let $v$ be a {\it derangement} of order~$m$ and
$w=x_1x_2\cdots x_n$ be a word from the shuffle class
$\Sh(0^{n-m}v)$ $(0\le n\le m)$, so that 
$v=x_{j_1}x_{j_2}\cdots x_{j_m}$ for 
$1\le j_1<j_2<\cdots <j_m\le n$. Let ``red" (``reduction") be the
increasing bijection of $\{j_1,j_2,\ldots,j_m\}$ onto the interval
$[\,m\,]$. Remember that a positive letter~$x_k$ of~$w$ is said to
be excedent (resp. subexcedent) if and only if
$x_k>\red k$ (resp. $x_k<\red k$). Accordingly, a letter is
non-subexcedent if it is either equal to~0 or excedent.

We define $n$ bijections $\phi_l$
$(1\le l\le n)$ of~$\Sh(0^{n-m}v)$ onto itself in the following
manner: for each~$l$ such that $n-m+1\le l\le n$ let
$\phi_l(w):=w$. When $1\le l\le n-m$, let $x_j$ denote the $l$-th
letter of~$w$, equal to~0, when~$w$ is read {\it from left to right}.
Three cases are next considered (by convention,
$x_0=x_{n+1}=+\infty$):

(1) $x_{j-1}$, $x_{j+1}$ both non-subexcedent;

(2) $x_{j-1}$ non-subexcedent, $x_{j+1}$ subexcedent; or
$x_{j-1}$, $x_{j+1}$ both subexcedent with $x_{j-1}>x_{j+1}$;

(3) $x_{j-1}$ subexcedent, $x_{j+1}$ non-subexcedent; or
$x_{j-1}$, $x_{j+1}$ both subexcedent with $x_{j-1}<x_{j+1}$.

\goodbreak
\medskip 
When case (1) holds, let $\phi_l(w):=w$.

When case (2) holds, determine the {\it greatest} integer
$k\ge j+1$ such that 
$$\eqalignno{
&\quad x_{j+1}<x_{j+2}<\cdots <x_k<\red(k),\cr
\noalign{\hbox{so that}}
\noalign{\vskip-5pt}
w&=x_1\cdots x_{j-1}\;0\; x_{j+1}\cdots
x_k\;x_{k+1}\cdots x_n\cr
\noalign{\hbox{and define:}}
\phi_l(w)&:=x_1\cdots x_{j-1}\;
x_{j+1}\cdots x_k\;0\;x_{k+1}\cdots x_n.\cr
}
$$

When case (3) holds, determine the {\it smallest}
integer $i\le j-1$ such that
$$\eqalignno{
&\quad\red(i)>x_i>x_{i+1}>\cdots >x_{j-1},\cr
\noalign{\hbox{so that}}
\noalign{\vskip-5pt}
w&=x_1\cdots x_{i-1}\ x_i\cdots x_{j-1}\;0\;
x_{j+1}\cdots x_n\cr
\noalign{\hbox{and define:}}
\phi_l(w)&:=x_1\cdots x_{i-1}\; 0\;
x_{i}\cdots x_{j-1}\;x_{j+1}\cdots x_n.\cr
}
$$

\goodbreak
It is important to note that $\phi_l$ has no action on the 0's
other than the $l$-th one. Then the mapping $\bfPhi$ in
Theorem~1.1 is defined to be the composition product
$$
\bfPhi:=\phi_1\phi_2\cdots \phi_{n-1}\phi_n.
$$

\smallskip
{\it Example}.\quad The following word $w$ has four zeros, so that
$\bfPhi(w)$ can be reached in four steps:
$$
\matrice{{\rm Id}&=&1\,&2\,&3\,&4\,&5\,&6\,&7\,&8\,&9\,&10&11&\cr
w&=&5&{\bf0}&1&2&{\bf0}&{\bf0}&3&6&{\bf0}&7&4\,&\  j=9, {\rm
apply}\
\phi_4, {\rm case\ (1)};\qquad\quad\cr
&&5&{\bf0}&1&2&{\bf0}&{\bf0}&3&6&{\bf0}&7&4\,&\  j=6, {\rm
apply}\
\phi_3, {\rm case\ (2)},\, k=7;\cr
&&5&{\bf0}&1&2&{\bf0}&3&{\bf0}&6&{\bf0}&7&4\,&\  j=5, {\rm
apply}\
\phi_2, {\rm case\ (3)},\, i=4;\cr
&&5&{\bf0}&1&{\bf0}&2&3&{\bf0}&6&{\bf0}&7&4\,&\  j=2, {\rm
apply}\
\phi_1, {\rm case\ (2)},\, k=3;\cr
\bfPhi(w)&=&5&1&{\bf0}&{\bf0}&2&3&{\bf0}&6&{\bf0}&7&4.&\cr
}$$
We have: $\RISE w=\Rise^\bullet \bfPhi(w)
=\{2,3,5,6,7,9,11\}$, as desired.

\medskip
To verify that $\bfPhi$ is bijective, we introduce a class of
bijections~$\psi_l$, whose definitions are parallel to the
definitions of the $\phi_l$'s. Let
$w=x_1x_2\cdots x_n\in \Sh(0^{n-m}v)$ $(0\le m\le n)$. For
each~$l$ such that $n-m+1\le l\le n$ let $\psi_l(w):=w$. When
$1\le l\le n-m$, let~$x_j$ denote the $l$-th letter of~$w$, equal
to~$0$, when~$w$ is read {\it from left to right}. Consider the
following three cases (remember that $x_0=x_{n+1}=+\infty$ by
convention):

$(1')=(1)$ $x_{j-1}$, $x_{j+1}$ both non-subexcedent;

$(2')$ $x_{j-1}$ subexcedent, $x_{j+1}$ non-subexcedent; or
$x_{j-1}$, $x_{j+1}$ both subexcedent with $x_{j-1}>x_{j+1}$;

$(3')$ $x_{j-1}$ non-subexcedent, $x_{j+1}$ subexcedent;
or $x_{j-1}$, $x_{j+1}$ both subexcedent with
$x_{j-1}<x_{j+1}$.

\goodbreak
\medskip
When case $(1')$ holds, let $\psi_l(w):=w$.

When case $(2')$ holds, determine the {\it smallest} integer
$i\le j-1$ such that
$$\eqalignno{
&\quad x_i<x_{i+1}<\cdots <x_{j-1}<\red(j-1),\cr
\noalign{\hbox{so that}}
w&=x_1\cdots x_{i-1}\ x_i\cdots x_{j-1}\;0\;
x_{j+1}\cdots x_n\cr
\noalign{\hbox{and define:}}
\psi_l(w)&:=x_1\cdots x_{i-1}\;0\;
x_{i}\cdots x_{j-1}\;x_{j+1}\cdots x_n.\cr
}
$$

When case $(3')$ holds, determine the {\it greatest} integer
$k\ge j+1$ such that 
$$\eqalignno{&\quad \red(j+1)>x_{j+1}>x_{j+2}>\cdots
>x_k,\cr
\noalign{\hbox{so that}}
w&=x_1\cdots x_{j-1}\;0\; x_{j+1}\cdots
x_k\;x_{k+1}\cdots x_n\cr
\noalign{\hbox{and define:}}
\psi_l(w)&:=x_1\cdots x_{j-1}\;
x_{j+1}\cdots x_k\;0\;x_{k+1}\cdots x_n.\cr
}
$$

\goodbreak
We now observe that when case (2)
(resp. (3)) holds for~$w$, then case~$(2')$ (resp.~$(3')$) holds for
$\phi_l(w)$. Also, when case $(2')$ (resp. $(3')$) holds for~$w$, then
case~(2) (resp. (3)) holds for $\psi_l(w)$. Therefore
$$
\phi_l\psi_l=\psi_l\phi_l=\hbox{Identity map}
$$
and the product 
$\bfPsi:=\psi_n\psi_{n-1}\cdots \psi_2\psi_1$ is the inverse
bijection of~$\Phi$.

\bigskip
\centerline{\bf 3. Verification of $\RISE w=\RISE^\bullet\bfPhi(w)$} 

\medskip
Let us introduce an alternate definition for~$\bfPhi$. Let $w$
belong to $\Sh(0^{n-m}v)$ and $w'$ be a nonempty left factor
of~$w$, of length~$n'$. Let~$w'$ have~$p'$ letters equal to~0. If
$p'\ge 1$, write
$$
w'=x_1\cdots x_{j-1}\;0^{h} x_{j+h}x_{j+h+1}\cdots
x_{n'},
$$
where $1\le h\le p'$, $x_{j-1}\not=0$ and where the right factor
$x_{j+h}x_{j+h+1}\cdots x_{n'}$ contains no~0. By
convention $x_{j+h}:=+\infty$ if $j+h=n'+1$. 
If $x_{j+h}$ is subexcedent
let~$k$ be the {\it greatest} integer
$k\ge j+h$ such that $x_{j+h}<x_{j+h+1}<\cdots <x_k<\red(k)$.
If $x_{j-1}$ is subexcedent 
let~$i$ be the {\it smallest}
integer
$i\le j-1$ such that $\red(i)>x_i>x_{i+1}>\cdots >x_{j-1}$.
Examine four cases:

(1) if $x_{j-1}$ and $x_{j+h}$ are both excedent,
let 
$$\eqalignno{
u:=x_1\cdots x_{j-1},\quad
u'&:=0^{h} x_{j+h}x_{j+h+1}\cdots x_{n'}\cr
\noalign{\hbox{and define}}
\theta(u')&:=u'.\cr}
$$

\goodbreak
(2) if $x_{j-1}$ is excedent and $x_{j+h}$ subexcedent, or if
$x_{j-1}$, $x_{j+h}$ are both subexcedent with $x_{j-1}>x_{j+h}$,
let
$$\eqalignno{\noalign{\vskip-2pt}
u:=x_1\cdots x_{j-1},\quad
u'&:=0^{h}\, x_{j+h}\cdots x_k\,x_{k+1}\cdots
x_{n'}\cr
\noalign{\hbox{and define}}
\noalign{\vskip-4pt}
\theta(u')&:=x_{j+h}\cdots x_k\,0^{h}\,
x_{k+1}\cdots
x_{n'}.\cr}
$$

(3) if $x_{j-1}$, $x_{j+h}$ are both subexcedent with
$x_{j-1}<x_{j+h}$,
let
$$
\eqalignno{u:=x_1\cdots x_{i-1},\quad
u'&:=x_i\cdots x_{j-1}0^{h} x_{j+h}\cdots x_k\,x_{k+1}\cdots
x_{n'}\cr
\noalign{\hbox{and define}}
\theta(u')&:=0\,x_i\cdots x_{j-1}x_{j+h}\cdots
x_k\,0^{h-1}\, x_{k+1}\cdots
x_{n'}.\cr}
$$

(4) if $x_{j-1}$ is subexcedent and 
$x_{j+h}$ excedent,
let
$$
\eqalignno{u:=x_1\cdots x_{i-1},\quad
u'&:=x_i\cdots x_{j-1}0^{h} x_{j+h}\cdots 
x_{n'}\cr
\noalign{\hbox{and define}}
\theta(u')&:=0\,x_i\cdots x_{j-1}\,0^{h-1}\, x_{j+h}\cdots
x_{n'}.\cr}
$$

\goodbreak
\noindent
By construction $w'=uu'$. Call it the {\it canonical factorization}
of~$w'$. In the three cases we evidently have:
$$\leqalignno{
\RISE u'&=\RISE^\bullet \theta(u').&(3.1)\cr
\noalign{\hbox{Define $\Theta(w')$ to be the
three-term sequence:}}
\Theta(w')&:=(u,u',\theta(u')).&(3.2)\cr
}
$$
The rightmost letter of~$u$ occurs as the $q$-th
letter of~$w'$, so that $q=j-1$ in cases~(1) and~(2) and $i-1$ in
cases~(3) and~(4). 

\goodbreak
\proclaim Lemma 3.1. We have
$$
\leqalignno{
\RISE x_q u'&=\Rise^\bullet x_q \theta(u'),&(3.3)\cr
\RISE x_q u'&=\Rise^\bullet 0\theta(u'), \quad\hbox{\rm if $x_q$ is subexcedent}.&(3.4)\cr
}
$$

\proof
Because of (3.1) we only have to study the two-letter factor $x_qx_{q+1}$.
First, let us prove identity (3.3).
There is nothing to do in case~(1). 
In case~(2) we have to verify
$\RISE x_{j-1}0=\RISE^\bullet x_{j-1}x_{j+h}$.
When $x_{j-1}$ is excedent and $x_{j+h}$ subexcedent, then
 $1\not\in \RISE x_{j-1}0$ and
$1\not\in \RISE^\bullet x_{j-1}x_{j+h}$. When $x_{j-1}$,
$x_{j+h}$ are both subexcedent with $x_{j-1}>x_{j+h}$, then
$1\not\in \RISE x_{j-1}0$ and
$1\not\in \RISE^\bullet x_{j-1}x_{j+h}$.

In cases (3) and (4) we have to verify
$\RISE x_{i-1}x_i=\RISE^\bullet x_{i-1}0$.
If $x_{i-1}=0$ (resp.
excedent), then $1\in \RISE x_{i-1}x_i$ (resp. 
$1\not\in\RISE x_{i-1}x_i$) and
$1\in \RISE^\bullet x_{i-1}0$ (resp.
$1\not\in \RISE^\bullet x_{i-1}0$).
When $x_{i-1}$ is subexcedent, then $x_{i-1}<x_i$ by definition
of~$i$. Hence $1\in \RISE x_{i-1}x_i$ and $1\in \RISE^\bullet
x_{i-1}0$.

\goodbreak
\smallskip
We next prove identity (3.4). In case (1) $z$ is always excedent,
so that identity (3.4) needs not to be considered.
In case~(2) we have to verify
$\RISE x_{j-1}0=\RISE^\bullet 0x_{j+h}$.
But if $x_{j-1}$ is subexcedent, then $x_{j+h}$ is also subexcedent, 
so that the above two sets are empty. 
In cases (3) and (4) we have to verify
$\RISE x_{i-1}x_i=\RISE^\bullet 00=\{1\}$.
But if $x_{i-1}$ is subexcedent, then $x_{i-1}<x_i$ by definition of $i$,
so that $1\in \RISE x_{i-1}x_i$.
\cqfd

\medskip
Now, if $w$ has~$p$ letters equal to~0 with $p\ge1$,
it may be expressed as the juxtaposition product
$$
w=w_1\,0^{h_1}\,w_2\,0^{h_2}\,\cdots\,
w_r\,0 ^{h_r}\,w_{r+1},
$$
where $h_1\ge 1$, $h_2\ge 2$, \dots~, $h_r\ge 1$ and where the
factors $w_1$, $w_2$, \dots~, $w_r$, $w_{r+1}$ contain no
0 and $w_2$, \dots~, $w_r$ are nonempty. We may define:
$\Theta(w):=(u_r,u'_r,\theta(u'_r))$, where $w=u_ru'_r$ is the
canonical factorization of~$w$.
 As~$u_r$ is a left factor of~$w$, we next define
$\Theta(u_r):=(u_{r-1},u'_{r-1},\theta(u'_{r-1}))$,
where $u_r=u_{r-1}u'_{r-1}$ is the canonical factorization
of~$u_r$, and successively
$\Theta(u_{r-1}):=(u_{r-2},u'_{r-2},\theta(u'_{r-2}))$
with $u_{r-1}=u_{r-2}u'_{r-2}$, \dots~,
$\Theta(u_2):=(u_1,u'_1,\theta(u'_1))$ with $u_2=u_1u'_1$, so that
$w=u_1u'_1u'_2\cdots u'_r$ and
$\bfPhi(w)=u_1\,\theta(u'_1)\theta(u'_2)\cdots \theta(u'_r)$.

It can be verified that 
$u_r\,\theta(u'_r)=\phi_{p-h_r+1}\cdots \phi_{p-1}\phi_p(w)$,
\hfil\break
$u_{r-1}\theta(u'_{r-1})\theta(u'_r)
=\phi_{p-h_r-h_{r-1}+1}\cdots \phi_{p-1}\phi_p(w)$, etc.
\medskip
With identities (3.1), (3.3) and (3.4) the proof of (1.5) is now completed.

\medskip
Again, consider the word~$w$ of the preceding example
$$
w=5\,0\,1\,2\,0\,0\,3\,6\,0\,7\,4,
$$
so that $r=3$, $h_1=1$, $h_2=2$, $h_3=1$, $w_1=5$,
$w_2=1\,2$ and $w_3=3\,6$, $w_4=7\,4$.
We have 
$$
\eqalignno{
\Theta(w)&=(u_3,u'_3,\theta(u'_3))
=(5\,0\,1\,2\,0\,0\,3\,6;\  0\,7\,4;\; 0\,7\,4);&\hbox{case (1)}\cr
\Theta(u_3)&=(u_2,u'_2,\theta(u'_2))
=(5\,0\,1;\  2\,0\,0\,3\,6;\  0\,2\,3\,0\,6);&\hbox{case (3)}\cr
\Theta(u_2)&=(u_1,u'_1,\theta(u'_1))
=(5;\ 0\,1;\ 1\,0);&\hbox{case (2)}\cr
\bfPhi(w)&=u_1\,\theta(u'_1)\theta(u'_2)\theta(u'_3)
=5\mid 1\,0\mid0\,2\,3\,
0\,6\mid 0\,7\,4.&\cr}
$$

\bigskip
\centerline{\bf 4. The transformation ${\bf F}_3$} 

\medskip
The bijection ${\bf F}_3$ we are now defining maps each shuffle
class $\Sh(0^{n-m}v)$ with $v$ an arbitrary word of length~$m$
$(0\le m\le n)$ onto itself. When $n=1$ the unique element of the
shuffle class is sent onto itself. Let $n\ge 2$ and assume that~${\bf
F}_3(w')$ has been defined for all words $w'$ with nonnegative
letters, of length $n'\le n-1$. Further assume that (1.6) and
(1.7) hold for all those words. Let
$w$ be a word of length~$n$. We may write
$$
w=w'a0^rb,
$$
where $a\ge 1$, $b\ge 0$ and $r\ge 0$. Three cases are being
considered:

(1) $a\le b$;\quad (2) $a>b$, $r\ge 1$;\quad
(3) $a>b$, $r=0$.

\medskip
In case (1) define:\quad ${\bf F}_3(w)={\bf F}_3(w'a0^rb):=
({\bf F}_3(w'a0^r))b$.

In case~(2) we may write ${\bf F}_3 (w'a0^r)=w''0$ by Property (1.7).
We then define
$$\leqalignno{
\gamma\,{\bf F}_3(w'a0^r)&:=0w'';\cr
{\bf F}_3(w)={\bf F}_3(w'a0^rb)&:=(\gamma\,{\bf F}_3(w'a0^r))b=0w''b.\cr}
$$
In short, add one letter~``0" to the left of ${\bf F}_3(w'a0^r)$, then
delete the rightmost  letter~``0"  and add $b$ to the right.

In case (3) remember that $r=0$. Write
$$
{\bf F}_3(w'a)=0^{m_1}x_1v_10^{m_2}x_2v_2\cdots 0^{m_k}x_kv_k,
$$
where $m_1\ge 0$, $m_2,\ldots,m_{k}$ are
all positive, then
$x_1$, $x_2$, \dots~, $x_k$ are positive letters and $v_1$,
$v_2$, \dots~, $v_k$ are words with positive letters, possibly
empty. Then define:
$$\leqalignno{
\delta\,{\bf F}_3(w'a)&:=
x_10^{m_1}v_1x_20^{m_2}v_2x_3\cdots x_k0^{m_k}v_k;\cr
{\bf F}_3(w)={\bf F}_3(w'ab)&:=(\delta\,{\bf F}_3(w'a))b.\cr}
$$
In short, move each positive letter occurring just after a 0-factor
of ${\bf F}_3(w'a)$ to the beginning of that 0-factor and add~$b$
to the right.
\goodbreak
\medskip

{\it Example}.

\vskip-33pt
$$
\eqalignno{
w&=0\,0\,0\,3\,1\,2\,2\,0\,0\,1\,3\cr
{\bf F}_3(0\,0\,0\,3)&=0\,0\,0\,3&\hbox{case (1)}\cr
{\bf F}_3(0\,0\,0\,3\,1)&
=\delta(0\,0\,0\,3)1=3\,0\,0\,0\,1&\hbox{case (3)}\cr
{\bf F}_3(0\,0\,0\,3\,1\,2\,2)&=3\,0\,0\,0\,1\,2\,2&\hbox{case (1)}\cr
{\bf F}_3(0\,0\,0\,3\,1\,2\,2\,0)&=
\delta(3\,0\,0\,0\,1\,2\,2)0=3\,1\,0\,0\,0\,2\,2\,0&\hbox{case (3)}\cr
{\bf F}_3(0\,0\,0\,3\,1\,2\,2\,0\,0)&
=\gamma(3\,1\,0\,0\,0\,2\,2\,0)0
=0\,3\,1\,0\,0\,0\,2\,2\,0&\hbox{case (2)}\cr
{\bf F}_3(0\,0\,0\,3\,1\,2\,2\,0\,0\,1)&
=\gamma(0\,3\,1\,0\,0\,0\,2\,2\,0)1
=0\,0\,3\,1\,0\,0\,0\,2\,2\,1\qquad&\hbox{case (2)}\cr
{\bf F}_3(0\,0\,0\,3\,1\,2\,2\,0\,0\,1\,3)&
=0\,0\,3\,1\,0\,0\,0\,2\,2\,1\,3.&\cr
}
$$
We have:
$\maj w=\maj (0\,0\,0\,3\,1\,2\,2\,0\,0\,1\,3)=4+7=11$ 
and
$\mafz  {\bf F}_3(w)=\mafz  (0\,0\,3\,1\,0\,0\,0\,2\,2\,1\,3)
=(1+2+5+6+7)-(1+2+3+4+5)+(1+4)=11$.

\medskip

By construction the rightmost letter is preserved by~${\bf F}_3$.
To prove (1.6) proceed by induction. Assume that
$\mafz w'a0^r=\mafz {\bf F}_3(w'a0^r)$ holds. In case~(1)
``maj" and ``mafz" remain invariant when~$b$ is
juxtaposed at the end. In case~(2) we have
$\maj w=\maj(w'a0^rb)=\maj(w'a0^r)$, but
$\mafz \gamma\,{\bf F}_3(w'a0^r)=
\mafz {\bf F}_3(w'a0^r)-|w'a0^r|_{\ge 1}$
and $\mafz\, (\gamma\,{\bf F}_3(w'a0^r))b
=\mafz \gamma\,{\bf F}_3(w'a0^r)+|w'a0^r|_{\ge 1}$,
where $|w'a0^r|_{\ge 1}$ denotes the number of {\it positive}
letters in $w'a0^r$. Hence (1.6) holds. In case~(3)
remember $r=0$. We have
$\maj (w'ab)=\maj( w'a)+|w'a|$, where $|w'a|$ denotes the length of
the word~$w'a$. But
$\mafz \delta {\bf F}_3  (w'a)=\mafz{\bf F}_3  (w'a)+ \zero 
(w'a)$ and $\mafz (\delta {\bf F}_3  (w'a))b
=\mafz \delta {\bf F}_3  (w'a)+|w'a|_{\ge 1}$.
The equality holds for $b=0$ and $b\geq 1$, as easily verified.
As $\zero (w'a)+|w'a|_{\ge 1}=|w'a|$, we have
$\mafz{\bf F}_3(w)= \mafz\, (\delta {\bf F}_3  (w'a))b=\mafz
{\bf F}_3  (w'a)+|w'a|=\maj w'ab=\maj w$. Thus (1.6) holds in the
three cases.

\medskip
To define the inverse bijection ${\bf F}_3^{-1}$ of ${\bf F}_3$
we first need the inverses~$\gamma^{-1}(w)$
and~$\delta^{-1}(w)$ for each word~$w$.
Let~$w=0w'$ be a word, whose first letter is~0.
Define~$\gamma^{-1}(w)$ to be the word derived from~$w$ by
deleting the first letter~0 and adding one
letter~``0" to the right of~$w$. Clearly,
$\gamma^{-1}\gamma=\gamma\gamma^{-1}$ is the identity
map.

Next, let $w$ be a word, whose first letter is positive.
Define~$\delta^{-1}(w)$ to be the word derived from~$w$
by moving each positive letter occurring just before a 0-factor
of~$w$ to the end of that 0-factor. Again
$\delta^{-1}\delta=\delta\delta^{-1}$ is the identity map.

We may write 
$$ 
w=cw'a0^rb,
$$
where $a\ge 1$, $b\ge 0$, $c\ge 0$ and $r\ge 0$. Three cases are being
considered:

(1) $a\le b$;\quad (2) $a>b$, $c= 0$;\quad
(3) $a>b$, $c\ge 1$.

In case (1) define:\quad ${\bf F}^{-1}_3(w):=
({\bf F}^{-1}_3(cw'a0^r))b$.

In case (2) define:\quad ${\bf F}^{-1}_3(w):=
(\gamma^{-1}({\bf F}^{-1}_3(cw'a0^r)))b$.

In case (3) define:\quad ${\bf F}^{-1}_3(w):=
(\delta^{-1}({\bf F}^{-1}_3(cw'a0^r)))b$.

\goodbreak
\medskip
We end this section by proving a property of the
transformation~${\bf F}_3$, which will be used in our next paper
{\tt[FoHa07]}.

\proclaim Proposition 4.1. Let $w$, $w''$ be two words with
nonnegative letters, of the same length. If $\Zero w=\Zero w''$
and $\DES\Pos w=\DES\Pos w''$, then
$\Zero{\bf F}_3(w)=\Zero{\bf F}_3(w'')$.

{\it Proof}.\quad
To derive ${\bf F}_3(w)$ (resp. ${\bf F}_3(w'')$)
from $w$ (resp. $w''$) we have to consider one of the 
three cases (1), (2) or (3), described above, at each step.
Because of the two conditions
$\Zero w=\Zero w''$ and $\DES\Pos w=\DES\Pos w''$, 
case $(i)$ ($i=1,2,3$) is used at the $j$-th step in the
calculation of ${\bf F}_3(w)$, if and only if the same case
is used at that $j$-th step for the calculation of
${\bf F}_3(w'')$.
Consequently the letters equal to $0$ are in the same places in both words
${\bf F}_3(w)$ and ${\bf F}_3(w'')$.
\cqfd

\medskip

By the very definition of
$\Phi:\sigma\buildrel \ZDer\over \mapsto w\buildrel \bfPhi\over
\mapsto w'\buildrel \ZDer^{-1}\over \mapsto \sigma'$, given in
(1.18) and of
$\rmF_3:\sigma\buildrel \ZDer\over \mapsto w\buildrel {\bf
F}_3\over
\mapsto w''\buildrel \ZDer^{-1}\over \mapsto \sigma''$, given in
(1.19) we have
$\Phi(\sigma)=\sigma$ and $\rmF_3(\sigma)=\sigma$ if~$\sigma$
is a derangement. In the next two tables we have calculated
$\Phi(\sigma)=\sigma'$ and $\rmF_3(\sigma)=\sigma''$ for the
fifteen non-derangement permutations~$\sigma$ of order~4.

\bigskip

{
\eightpoint
$$
\vbox{\offinterlineskip\halign{
\vrule#\hfil\strut\ &\hfil$#$\hfil\ \vrule&\ \hfil$#$\hfil\ \vrule
&\ \hfil$#$\hfil\ \vrule&\ \hfil$#$\hfil\ \vrule&\ \hfil$#$\hfil\
\vrule &\ \hfil$#$\hfil\ \vrule&\ \hfil$#$\hfil\ \vrule&
\ \hfil$#$\hfil\ \vrule&\ \hfil$#$\hfil\ \vrule&\ \hfil$#$\hfil &\
\hfil#\vrule\cr
\noalign{\hrule}
&\fix\sigma&\Der\sigma&\hbox{\sevenrm RIZE}\,\sigma&
\sigma&w
&w'&\sigma'
&\hbox{\sevenrm
RISE}\,\sigma'&\Der \sigma'&\fix \sigma'&\cr
\noalign{\hrule}
&4&e&1,2,3,4&1234&{\bf0}\,{\bf0}\,{\bf0}\,{\bf0}\,
&{\bf0}\,{\bf0}\,{\bf0}\,{\bf0}\,&1234
&1,2,3,4&e&4&\cr
\noalign{\hrule}
&&&1,2,4&1243&{\bf0}\,{\bf0}\,
2\,1&{\bf0}\,{\bf0}\,2\,1&1243&1,2,4&&&\cr
&&&1,4&1324&{\bf0}\,2\,1\,{\bf0}&{\bf0}\,2{\bf0}\,1
&1432&1,4&&&\cr 
&&&1,3,4&1432&{\bf0}\,2\,{\bf0}\,1&{\bf0}\,2\,1\,{\bf0}
&1324&1,3,4&&&\cr
&2&21&3,4&2134&2\,1\,{\bf0}\,{\bf0}&2\,{\bf0}\,1\,{\bf0}
&3214&3,4&21&2&\cr
&&&2,4&3214&2\,{\bf0}\,1\,{\bf0}&2\,{\bf0}\,{\bf0}\,1
&4231&2,4&&&\cr
&&&2,3,4&4231&2\,{\bf0}\,{\bf0}\,1&2\,1\,{\bf0}\,{\bf0}
&2134&2,3,4&&&\cr
\noalign{\hrule}
&&&1,2,4&1342&{\bf0}\,2\,3\,1\,&{\bf0}\,2\,3\,1\,
&1342&1,2,4&&&\cr
&&&1,4&2314&2\,3\,1\,{\bf0}&2\,3\,{\bf0}\,1
&2431&1,4&&&\cr
&1&231&1,3,4&2431&2\,3\,{\bf0}\,1&2\,3\,1\,{\bf0}
&2314&1,3,4&231&1&\cr
&&&2,4&3241&2\,{\bf0}\,3\,1&2\,{\bf0}\,3\,1
&3241&2,4&&&\cr
\noalign{\hrule}
&&&1,3,4&1423&{\bf0}\,3\,1\,2&{\bf0}\,3\,1\,2
&1423&1,3,4&&&\cr
&1&312&2,4&3124&3\,1\,2\,{\bf0}&3\,1\,{\bf0}\,2
&4132&2,4&312&1&\cr
&&&3,4&4132&3\,1\,{\bf0}\,2&3\,{\bf0}\,1\,2
&4213&3,4&&&\cr
&&&2,3,4&4213&3\,{\bf0}\,1\,2&3\,1\,2\,{\bf0}
&3124&2,3,4&&&\cr
\noalign{\hrule}
}}
$$

\centerline{Calculation of $\sigma'=\Phi(\sigma)$}
}

{\eightpoint

$$
\vbox{\offinterlineskip\halign{
\vrule#\hfil\strut\ &\hfil$#$\hfil\ \vrule&\ \hfil$#$\hfil\ \vrule
&\ \hfil$#$\hfil\ \vrule&\ \hfil$#$\hfil\ \vrule&\ \hfil$#$\hfil\
\vrule &\ \hfil$#$\hfil\ \vrule&\ \hfil$#$\hfil\ \vrule&
\ \hfil$#$\hfil\ \vrule&\ \hfil$#$\hfil\ \vrule&\ \hfil$#$\hfil &\
\hfil#\vrule\cr
\noalign{\hrule}
&\fix\sigma&\Der\sigma&\maz\sigma&
\sigma&w
&w''&\sigma''
&\maf\sigma''&\Der \sigma''&\fix \sigma''&\cr
\noalign{\hrule}
&4&e&0&1234&{\bf0}\,{\bf0}\,{\bf0}\,{\bf0}\,
&{\bf0}\,{\bf0}\,{\bf0}\,{\bf0}\,&1234
&0&e&4&\cr
\noalign{\hrule}
&&&3&1243&{\bf0}\,{\bf0}\,
2\,1&2\,{\bf0}\,{\bf0}\,1&4231&3&&&\cr
&&&5&1324&{\bf0}\,2\,1\,{\bf0}&2\,1\,{\bf0}\,{\bf0}
&2134&5&&&\cr 
&&&2&1432&{\bf0}\,2\,{\bf0}\,1&{\bf0}\,2\,{\bf0}\,1
&1432&2&&&\cr
&2&21&3&2134&2\,1\,{\bf0}\,{\bf0}&{\bf0}\,2\,1\,{\bf0}
&1324&3&21&2&\cr
&&&4&3214&2\,{\bf0}\,1\,{\bf0}&2\,{\bf0}\,1\,{\bf0}
&3214&4&&&\cr
&&&1&4231&2\,{\bf0}\,{\bf0}\,1&{\bf0}\,{\bf0}\,2\,1
&1243&1&&&\cr
\noalign{\hrule}
&&&3&1342&{\bf0}\,2\,3\,1\,&2\,{\bf0}\,3\,1\,
&3241&3&&&\cr
&&&5&2314&2\,3\,1\,{\bf0}&2\,3\,1\,{\bf0}
&2314&5&&&\cr
&1&231&2&2431&2\,3\,{\bf0}\,1&{\bf0}\,2\,3\,1
&1342&2&231&1&\cr
&&&4&3241&2\,{\bf0}\,3\,1&2\,3\,{\bf0}\,1
&2431&4&&&\cr
\noalign{\hrule}
&&&2&1423&{\bf0}\,3\,1\,2&3\,{\bf0}\,1\,2
&4213&2&&&\cr
&1&312&4&3124&3\,1\,2\,{\bf0}&3\,1\,2\,{\bf0}
&3124&4&312&1&\cr
&&&3&4132&3\,1\,{\bf0}\,2&3\,1\,{\bf0}\,2
&4213&3&&&\cr
&&&1&4213&3\,{\bf0}\,1\,2&{\bf0}\,3\,1\,2
&3124&1&&&\cr
\noalign{\hrule}
}}
$$

\centerline{Calculation of $\sigma''=\rmF_3(\sigma)$}
}

\vfill\eject
\vglue 2mm


{
\eightpoint
\def\thevskip{\smallskip}
\centerline{\bf References} 
\bigskip

\article CHZ97|R. J. Clarke, G.-N. Han, J. Zeng|A combinatorial
interpretation of the Seidel generation of $q$-derangement numbers|%
Annals of Combinatorics|4|1997|313--327|
\thevskip

\divers FoHa07|Dominique Foata, Guo-Niu Han|Fix-Mahonian
Calculus, II: further statistics, preprint, 13 p., {\oldstyle 2007}|
\thevskip

\livre GaRa90|George Gasper,
Mizan Rahman|Basic Hypergeometric Series|London,
Cambridge Univ. Press, {\oldstyle 1990}  ({\sl Encyclopedia of
Math. and Its Appl.}, {\bf 35})|
\thevskip

\article GR93|I. Gessel, C. Reutenauer|Counting Permutations with
Given Cycle Structure and Descent Set| J. Combin. Theory Ser.
A|64|1993|189--215|

\livre Lo83|M. Lothaire|Combinatorics on Words|Addison-Wesley,
London {\oldstyle 1983} (Encyclopedia of Math. and its Appl., {\bf
17})|
\thevskip

\divers ShWa06|John Shareshian, Michelle L. Wachs|$q$-Eulerian
Polynomials: Excedance Number and Major Index, preprint, 14 p.,
{\oldstyle 2006}|
\thevskip

}

\bigskip\bigskip
\hbox{\vtop{\halign{#\hfil\cr
Dominique Foata \cr
Institut Lothaire\cr
1, rue Murner\cr
F-67000 Strasbourg, France\cr
\noalign{\smallskip}
{\tt foata@math.u-strasbg.fr}\cr}}
\qquad
\vtop{\halign{#\hfil\cr
Guo-Niu Han\cr
I.R.M.A. UMR 7501\cr
Universit\'e Louis Pasteur et CNRS\cr
7, rue Ren\'e-Descartes\cr
F-67084 Strasbourg, France\cr
\noalign{\smallskip}
{\tt guoniu@math.u-strasbg.fr}\cr}}}

\vfill\eject

\bye

First, if
$w=y_1y_2\cdots y_n$ is a word, whose letters are integers, having
no repetitions. Say that~$w$ is a {\it desarrangement} if
$y_1>y_2>\cdots >y_{2k}$ and
$y_{2k}<y_{2k+1}$ for some~$k\ge 1$. By convention,
$y_{n+1}=\infty$. We could also say that the {\it leftmost
trough} of~$w$ occurs at an {\it even} position. This notion
was introduced by D\'esarm\'enien {\tt[De84]} and elegantly used
in a subsequent paper {\tt[DeWa88]}. A further refinement is due
to Gessel {\tt[Ge91]}. Each desarrangement which is also a
permutation from ${\goth S}_n$ is called a {\it desarrangement} of
order~$n$. Let~$K_n$ be the set of all those desarrangements. If
$\sigma$ is in $K_n$, the inverse~$\sigma^{-1}$ of~$\sigma$
is said to be an {\it idesarrangement} of order~$n$. The set of all
those idesarrangements is denoted by $I\kern-2pt K_n$. Clearly, an
idesarrangement of order~$n$ is a permutation
$\sigma=\sigma(1)\sigma(2)\cdots \sigma(n)$ such that 
$(2k)(2k-1)\cdots 2\,1$, but {\it not} $(2k+1)(2k)(2k-1)\cdots 2\,1$,
is a subword of~$\sigma$ for some $k\ge 1$.

Let $\sigma=\sigma(1)\sigma(2)\cdots\sigma(n)$ be a permutation.
Unless~$\sigma$ is increasing, there is always a nonempty right
factor of~$\sigma$ which is a desarrangement. It then makes sense
to define~$\sigma^d$ as the {\it longest} such a right factor.
Hence,~$\sigma$ admits a unique factorization
$\sigma=\sigma^p \sigma^d$, called the {\it pixed
factorization}, where $\sigma^p$ is  {\it increasing}
and~$\sigma^d$ is the longest right factor of~$\sigma$ which is
a desarrangement. The set (resp. number) of the letters
in~$\sigma^p$ is denoted by $\PIX\sigma$ (resp. $\pix\sigma$). 

Now write the inverse $\sigma^{-1}$ of~$\sigma$ as a word
$\sigma^{-1}(1)\sigma^{-1}(2)\cdots\sigma^{-1}(n)$. We have:
$\sigma^{-1}(i)\ge \pix+1$ if and only if~$i\in [n]\setminus
\PIX\sigma$. It makes sense to define the word
$Z^{\Idesar}(\sigma)=x_1x_2\cdots x_n$ by
$$
x_i:=\cases{0,&if $i\in \PIX\sigma$;\cr
\sigma^{-1}(i)-\pix\sigma,&if $i\in [n]\setminus
\PIX\sigma$.\cr}
$$
The subword of $Z^{\Idesar}(\sigma)=w$ obtained by deletion of the
zeros is an idesarrangement of order~$n-\pix\sigma$, as is readily
seen, which is further denoted by $\Idesar\sigma$.

\medskip
{\it Example}.\quad
Let $\sigma=2\,4\,8\,5\,3\,9\,1\,7\,6$ be a permutation. Its pixed
factorization is $\sigma^p\sigma^d$ with $\sigma^p=2\,4\,8$ and
$\sigma^d=5\,3\,9\,1\,7\,6$. Thus, $\PIX\sigma=\{2,4,8\}$,
$\pix\sigma=3$. Furthermore,
$\sigma^{-1}=7\,1\,5\,2\,4\,9\,8\,3\,6$,
$Z^{\Idesar}(\sigma)=4\,0\,2\,0\,1\,6\,5\,0\,3$
and $\Idesar\sigma=4\,2\,1\,6\,5\,3$.

If $\sigma^d=\sigma(n-m+1)\sigma(n-m+2)\cdots
\sigma(n)$ and if ``$\red$" is the increasing bijection mapping
the set $\{\sigma(n-m+1),\sigma(n-m+2),\ldots,
\sigma(n)\}$ onto $\{1,2,\ldots,m\}$, define
$$\leqalignno{
\Desar\sigma&:=\red\sigma(n-m+1)\,\red\sigma(n-m+2)\,\ldots\,
\red\sigma(n);&(1.9)\cr
\mag\sigma&:=\sum_{i\in \hbox{\sevenrm PIX}(\sigma)}
\kern-5pt i-\sum\limits_{i=1}^{\pix\sigma}i+\imaj\circ
\Desar\sigma.&(1.10)\cr}
$$

Note that $\Desar\sigma$ is a desarrangement and belongs to
${\goth S}_m$, for short, a desarrangement of order~$m$.
Also note that $\sigma$ is fully characterized by the 
pair $(\PIX\sigma,\Desar\sigma)$, which will called the {\it pixed
decomposition} of~$\sigma$.


\bye